\documentclass[11pt]{article}
\usepackage{amsmath,amsfonts,amssymb,mathrsfs,cases,mathdots,colordvi,url}

\newtheorem{lemma}{Lemma}[section]
\newtheorem{thm}{Theorem}[section]
\newtheorem{proposition}{Proposition}[section]

\newtheorem{example}{Example}[section]

\renewcommand{\Box}{\rule{2.2mm}{2.2mm}}
\newcommand{\BOX}{\hfill \Box}
\newcommand{\K}{{\cal K}}

\newcommand{\ds}{\displaystyle}
\date{December 13, 2015}

\begin{document}

\title{Exact formulas for the proximal/regular/limiting normal cone of the second-order cone complementarity set}
\author{Jane J. Ye\thanks{Department of Mathematics
and Statistics, University of Victoria, Victoria, B.C., Canada V8W 2Y2, e-mail: janeye@uvic.ca. The research of this author was partially
supported by NSERC.}\  \ \ \ and \ \ Jinchuan Zhou\thanks{Department of Mathematics, School of Science, Shandong University of Technology,
 Zibo 255049, P.R. China, e-mail: jinchuanzhou@163.com. This author's work is supported by National
 Natural Science Foundation of China (11101248, 11271233).}}

\maketitle
\begin{abstract}

 The proximal, regular and limiting  normal cones to the second-order cone complementarity set play   important roles in studying mathematical programs with second-order cone complementarity constraints, second-order cone programs, and the  second-order cone complementarity problems. It is needed in the first-order optimality conditions for mathematical programs with second-order cone complementarity constraint, the second-order subdifferential criteria in characterizing    the full stability for  second-order cone programs and  second-order cone complementarity problems, as well as in the characterizing the pseudo-Lipschitz continuity  of the solution mapping to  parametric second-order cone complementarity problems.  In this paper we establish explicit formulas for the proximal, regular, and limiting normal cone of the second-order cone complementarity set.

\vskip 10 true pt

\vskip 10 true pt


\end{abstract}

\section{Introduction}

 Let $X$ be a finite dimensional space and $\Theta\subset X$ be a convex set. We call
 $$\Omega:=\{(x,y)|\ x\in \Theta,\ y\in \Theta,\ \langle x,y\rangle =0\}$$
 a complementarity set associated with $\Theta$ or simply a complementarity set.
 Note that $\Omega$ is a cone whenever $\Theta$ is cone and in this case we may also call  $\Omega$ a complementarity cone.  Due to the existence of the complementarity condition, a complementarity set is always nonconvex and hence  is a difficult subject to study in  the variational analysis. Compared with results  for convex cones such as the second-order cone and the semidefinite matrix cone, so far there is not much research done in variational analysis for the complementarity set yet.

 \medskip

 Normal cones of the complementarity set play important roles in optimality conditions and stability analysis of optimization and equilibrium problems. For example, an optimization problem where some of the constraints  are in  the form of
the complementarity system
\begin{equation}
\Theta \ni G(z)\perp H(z) \in \Theta \label{equilibrium}
\end{equation}
 does not satisfy the classical constraint qualifications (see e.g. \cite{DSY14,yzz}).  To deal with this difficulty, one can reformulate
 (\ref{equilibrium})
as
\begin{equation} (G(z),H(z)) \in \Omega \label{set}
\end{equation} since as far as constraint qualifications concerned,  a constraint  in the form of (\ref{set})  is much easier to deal with than the original constraint in the form of (\ref{equilibrium}).
Based on this reformulation the stationary condition involving  the limiting normal cone and the proximal/regular normal cone of $\Omega$ is referred to as an M-and S-stationary condition  respectively (e.g. \cite{DSY14,y1999,ye05}).
It is well-known that the stability of a minimizer of a second-order smooth function is strongly associated with the positive-definiteness of the Hessian matrix of the function. Using the indicator function $\delta_\Theta$, a constrained optimization problem
\begin{eqnarray}
\min \ f(z) \ \ s.t. \ \ g(z) \in \Theta \label{constrained}
\end{eqnarray}
can be considered as a unconstrained optimization problem:
\begin{equation} \min \  f(z)+\delta_\Theta(g(z)).\label{indicator}  \end{equation}
The unconstrained optimization (\ref{indicator}), however, has an extended-valued objective function. Recent progresses in variational analysis show that stability of the problem (\ref{constrained}) can be characterized by using the second-order subdifferential of the objective function in (\ref{indicator});
see \cite[Theorem 5.6]{MNR}.  To calculate the generalized Hessian/the second-order subdifferential of the objective function, one needs to calculate the second-order subdifferential of the indicator function $\delta_\Theta$.
Given an element $\bar{y}$ lying in the limiting subdifferential of the indicator function $\partial \delta_\Theta (\bar{x})$, the second-order subdifferential of $\delta_\Theta $ is the set-valued mapping $\partial^2 \delta_\Theta(\bar x,\bar y)(\cdot)$
defined by
$$\partial^2 \delta_\Theta (\bar x,\bar y) (y^*):=\{x^*|(x^*,-y^*)\in N_{gph N_\Theta}(\bar x,\bar y)\} \quad \mbox{ for all } y^*\in X,$$
through the limiting normal cone of the graph of the limiting normal cone $N_\Theta$.
If $\Theta$ is a self-dual cone, then
$$
 (x,y)\in gph N_\Theta \Longleftrightarrow (x,-y)\in \Omega.
$$
Hence calculating the second-order subdifferential of $\delta_\Theta $ can be done by calculating the normal cone to the complementarity set:
$$N_{gph N_\Theta}(x,y)=\begin{bmatrix}
I & 0 \\
0 & -I
\end{bmatrix}N_\Omega (x,-y),$$
where $I$ the identity matrix of appropriate size.
Moreover using the second-order subdifferential of the indicator function $\delta_\Theta$, one can characterize the pseudo-Lipschitz continuity of the solution mapping to the complementarity systems (\ref{equilibrium}) with parameter $p$ in the form
$$S( p):=\{z|\  0\in H(p,z)+\partial \delta_\Theta(G(p,z))\};$$
see \cite[Theorem 5.1]{BorisOutrata}.

\medskip

 Some results have been given for $\Omega$ when $\Theta$ is a special convex cone. For example, i) in the case where $X=\mathbb{R}^n$ and $\Theta=\mathbb{R}_+^n$, the proximal normal cone and the limiting normal cone formula are well-known; see  \cite[Proposition 2.7]{y1999} and \cite[Proposition 3.7]{y2000} respectively.  Moreover it is  easy to show that the proximal normal cone coincides with the regular normal cone. ii) In the case where $X={\cal S}^n$ and $\Theta={\cal S}_+^n$, the positive semidefinite matrix cone, the proximal normal cone and the limiting normal cone formula are given in \cite[Proposition 3.2]{DSY14} and \cite[Theorem 3.1]{DSY14}  respectively.  Moreover it was  shown that the proximal normal cone coincides with the regular normal cone \cite[Page 551]{DSY14}.


\medskip

In this paper we derive exact formulas for proximal/regular and limiting normal cone for  the complementarity set in the case where $\Theta $ is equal to $\cal K$, the $m$-dimensional second-order cone defined by
\[
 {\cal K}:=\{x=(x_1,x_2)\in \mathbb{R}\times \mathbb{R}^{m-1}|\ x_1\geq \|x_2\|\}
 \]
  where $\|\cdot\|$ denotes the Euclidean norm. Such formulas will be useful to study the optimality conditions for mathematical programs with second-cone complementarity constraints and stability analysis of the second-order cone programming \cite{MOS}.
 By the definition of the metric projection operator $\Pi_{\cal K}$, it is easy to see that
\begin{equation}
(x,y)\in \Omega \Longleftrightarrow x=\Pi_{\cal K}(x-y).\label{relation-1}
\end{equation}
In \cite{OS08},
Outrata and Sun  derived the formulas for  the directional derivatives, the regular and the limiting coderivatives  of the metric projection.  Based on these formulas,  Liang, Zhu and Lin \cite{LZL14} tried to derive exact expressions for the regular and the limiting normal cones of the second-order cone complementarity set.
Unfortunately, there are some gaps in their expressions of the regular and the limiting normal cones.
In this paper we fill in these gaps by deriving the correct exact expressions for the regular and limiting normal cone of the second-order cone complementary set.
In addition, we further study the proximal normal cone and show that
the regular and the proximal normal cones coincide with each other.


 \medskip

\section{Preliminaries}
 In this section we summarize some background materials on variational analysis and second-order cone which will be used in the following analysis. Detailed discussions on these subjects can be found in \cite{AG03,c,clsw,LVBL98,m1,rw}.

 \medskip

  Let $C\subset \mathbb{R}^n$.
   $x^* \in {\rm cl}C$, the proximal normal cone and the regular/Fr\'{e}chet normal cone of $C$ at ${x^*}$ are defined as
  \begin{eqnarray*}
 N^\pi_C ( x^*)&:=& \{v\in \mathbb{R}^n |\  \exists\, M>0\, \mbox{ such that } \langle v, x-x^*\rangle \leq M\|x-x^*\|^2 \ \ \forall \, x\in C\} \label{eq:def-proximal-norm}\\
 \widehat{{N}}_C(x^*) &:=& \{v\in \mathbb{R}^n\,|\
 \langle v, x-x^*\rangle \leq o(\|x-x^*\|) \ \forall  x\in C\}\label{normals}
 \end{eqnarray*}
 respectively.
 The limiting/Mordukhovich normal cone is defined as
\begin{equation*}\label{eq:def-limit-norm}
N_C (x^*):= \{\lim_{i\rightarrow \infty} \zeta_i|\  \zeta_i \in N^\pi_C(x_i),\ \   x_i \rightarrow x^*,\ \  x_i\in C\}= \{\lim_{i\rightarrow \infty} \zeta_i|\  \zeta_i \in \widehat{N}_C(x_i),\ \   x_i \rightarrow x^*\ \  x_i\in C\}.
\end{equation*}
 Let $\Phi: \mathbb{R}^n \rightrightarrows \mathbb{R}^m$ be a set-valued map and $(x^*,y^*) \in {\rm gph } \Phi$, where ${\rm gph } \Phi$ denotes the graph of $\Phi$. The regular coderivative and the limiting (Mordukhovich) coderivative of $\Phi$ at $(x^*,y^*)$ are the set-valued maps defined by
\begin{eqnarray*}
\widehat{D}^*\Phi(x^*,y^*)(v)&:=&\{ u\in \mathbb{R}^n| (u,-v) \in \widehat{N}_{{\rm gph} \Phi}(x^*,y^*)\},
\\
{D}^*\Phi(x^*,y^*)(v)&:=&\{ u\in \mathbb{R}^n| (u,-v) \in {N}_{{\rm gph} \Phi}(x^*,y^*)\}
\end{eqnarray*}
respectively. We omit $y^*$ in the coderivative notation if the set-valued map $\Phi$  is single-valued at ${x^*}$.
Moreover if $\Phi$ is a continuously differentiable single-valued map, then
 \[
 \widehat{D}^*\Phi(x^*)={D}^*\Phi(x^*)={\cal J}\Phi(x^*),
\]
where ${\cal J}\Phi(x^*)$ denotes the Jacobian matrix of $\Phi$ at $x^*$.

\medskip


 The topological interior and the boundary of ${\cal K}$ are
\begin{eqnarray*}
{\rm int} {\cal K}= \{(x_1,x_2) \in \mathbb{R} \times \mathbb{R}^{m-1}|x_1>\|x_2\|\} \ \ {\rm and}\ \  {\rm bd} {\cal K}=\{(x_1,x_2) \in \mathbb{R} \times \mathbb{R}^{m-1}|x_1=\|x_2\|\},
\end{eqnarray*}
respectively. For any given vector $x:=(x_1,x_2)\in \mathbb{R}\times \mathbb{R}^{m-1}$, it can be decomposed as
$$x=\lambda_1(x)c_1(x)+\lambda_2(x)c_2(x),$$
where $\lambda_i(x)$ and $c_i(x)$ for $i=1,2$ are the spectral values and the associated spectral vectors of $x$ given by
\[
 \lambda_i(x)=x_1+(-1)^i\|x_2\|  \quad  {\rm and}\quad
 c_i(x)=\left\{\begin{array}{ll}\frac{1}{2}(1, (-1)^i\bar{x}_2) & {\rm if} \ \ x_2\neq 0\\
\frac{1}{2}(1,w) & {\rm if}\ \ x_2=0\end{array}\right.
\]
with $\bar{x}_2:={x_2}/\|x_2\|$ and $w$ being any vector in $\mathbb{R}^{m-1}$ satisfying $\|w\|=1$.
For $x\in \mathbb{R}^m$, let $\Pi_{\cal K}(x)$ be the metric projection  of $x$ onto ${\cal K}$. Then by \cite{FLT},  it can be calculated by
\begin{equation}
\Pi_{\cal K}(x)=(\lambda_1(x))_+c_1(x)+(\lambda_2(x))_+c_2(x).\label{proj}
\end{equation}

  As we will  show  in the following proposition,  the expressions of the regular and the limiting normal cone for the complementarity set can be derived from  the expression for the coderivatives of the metric projection operator.
\begin{proposition}\label{Prop2.2} Let $(x,y)\in \Omega:=\{(x,y)|x\in {\cal K}, y\in {\cal K}, x^Ty=0\}$. Then
\begin{eqnarray}
 \widehat{N}_\Omega (x,y)&=&\bigg\{(u,v)|\ -v\in \widehat{D}^*\Pi_{\cal K}(x-y)(-u-v)\bigg\}\label{regular-1},\\
 {N}_\Omega (x,y)&=&\bigg\{(u,v)|\ -v\in {D}^*\Pi_{\cal K}(x-y)(-u-v)\bigg\}\label{regular-2} .
 \end{eqnarray}
 \end{proposition}
 \noindent {\bf Proof.} By (\ref{relation-1}), $\Omega$ can be rewritten as
$ \Omega
 = \{(x,y)|\ (x-y,x)\in {\rm gph} \Pi_{\cal K}\}.$
The desired results follows from  applying  the change of coordinate formula in \cite[Exercise 6.7]{rw}.
\BOX

\medskip
Finally, we recall other notations that will be used throughout the paper. The inner product of two vectors $x,y$ is denoted by $x^Ty$ or $\langle x,y \rangle$. For any $t\in \mathbb{R}$, define $t_+:=\max\{0,t\}$ and $t_-:=\min\{0,t\}$. For $x=(x_1,x_2) \in \mathbb{R}\times \mathbb{R}^{m-1}$, we write its reflection about the $x_1$ axis as  $\hat{x}:=(x_1,-x_2)$. Given a vector $x$, denote by $\mathbb{R}x$ the set $\{tx|\ t\in \mathbb{R}\}$. $\mathbb{R}_+x$ and $\mathbb{R}_{++}x$ where $\mathbb{R}_+:=[0,\infty)$ and $\mathbb{R}_{++}:=(0,\infty)$ are similarly defined. The polar cone of a vector $v$ is $v^\circ:=\{x|\ x^Tv\leq 0\}$. For a differentiable  mapping $H:\mathbb{R}^n\to \mathbb{R}^m$ and a vector $x\in \mathbb{R}^n$, we denote by ${\cal J}H(x)$ the Jacobian matrix of $H$ at $x$ and $\nabla H(x):={\cal J}H(x)^T$. For a single-valued Lipschitz continuous map $\Phi:\mathbb{R}^n \rightarrow \mathbb{R}^m$, we denote  B(ouligand)-subdifferential   by
 $\partial_B \Phi( x)$ and $\Phi'(x;h)$ the directional derivative of $\Phi$ at $x$ in direction $h$.

\section{Expression of the regular normal cone}

 In \cite[Proposition 2.2]{LZL14},  Liang, Zhu and Lin gave a formula for the regular and  limiting  normal cone of $\Omega$. Their formula for the case where $(x,y)\in \Omega$ with $x,y\in {\rm bd}{\cal K}\backslash \{0\}$ is the following:
\begin{equation}
 \widehat{N}_\Omega(x,y)=\{(u,v)\in \mathbb{R}^m\times \mathbb{R}^m|\ u\in \mathbb{R} \hat{x}, v\in \mathbb{R}\hat{y} \}.
 \label{we1}
 \end{equation}
 The following example shows that  formula (\ref{we1}) is incorrect when the dimension $m$ is greater than $2$. In the meantime, the example illustrates our new formula.
 \begin{example}
 Take $x=(1,1/\sqrt{2},1/\sqrt{2})$ and $y=(2,-\sqrt{2},-\sqrt{2})$. It is easy to see that $(x,y)\in \Omega$ with $x,y\in {\rm bd}{\cal K}\backslash \{0\}$, and $y=2\hat{x}$. Let $u=(1/\sqrt{2},-1,0)$
 and $v=(1/(2\sqrt{2}),0,1/2)$. Since $x-y\in \mathbb{R}^3\backslash (-{\cal K}\cup {\cal K})$, by \cite[Lemma 1(i) and Theorem 1(i)]{OS08}, we have
 \[
{\cal J}\Pi_{\cal K} (x-y)
=\begin{bmatrix}
 \frac{1}{2} & \frac{1}{2\sqrt{2}} & \frac{1}{2\sqrt{2}} \\
 \frac{1}{2\sqrt{2}} & \frac{5}{12} & \frac{1}{12} \\
 \frac{1}{2\sqrt{2}} & \frac{1}{12} & \frac{5}{12}
 \end{bmatrix} \]
 and hence
$
 \widehat{D}^*\Pi_{\cal K}(x-y)(-u-v)
 =
 (-\frac{1}{2\sqrt{2}}, 0 , -\frac{1}{2})=-v.
$
By Proposition \ref{Prop2.2}, it follows that  $(u,v)\in \widehat{N}_\Omega (x,y)$. However since  $u\notin \mathbb{R} \hat{x}$ and $v\notin \mathbb{R} \hat{y}$, formula (\ref{we1}) is incorrect. In fact, according to our formula to be derived in Theorem \ref{formula-regular normal conenew}, $(u,v)$ is an element of the regular normal cone since $u\perp x, v\perp y$ and $x_1\hat{u}+y_1v=\sqrt{2}(1,\frac{1}{\sqrt{2}}, \frac{1}{\sqrt{2}})\in \mathbb{R}x$.
 \end{example}
In the following result, we revise the formula for the regular normal cone obtained in \cite[Proposition 2.2]{LZL14} for  the case where
$x,y \in {\rm bd}{\cal K}\backslash\{0\} , x^Ty=0$.
It is easy to see that  when $m=2$, the condition $ u\perp x , \ v\perp y,  \ x_1\hat{u}+y_1v\in \mathbb{R} x$ is equivalent to $ u\perp x , \ v\perp y$, which in turn is equivalent to $ u\in \mathbb{R} \hat{x}, v\in \mathbb{R}\hat{y} $. Hence when $m\leq 2$, our regular normal cone formula is the same as the one given in \cite[Proposition 2.2]{LZL14}.
 \begin{thm}\label{formula-regular normal conenew}
 Let $(x,y)\in \Omega:=\{(x,y)|x\in {\cal K}, y\in {\cal K}, x^Ty=0\}$. Then
 \begin{eqnarray*}
 \widehat{N}_\Omega(x,y)=&   \left  \{ \begin{array}{ll}
\{(u,v)| u\in \mathbb{R}^m, \ v=0 \} & \ {\rm if} \ x=0,\ y\in {\rm int}{\cal K}; \\
\{(u,v)| u=0, v\in \mathbb{R}^m  \}  & \ {\rm if}\ x\in {\rm int}{\cal K} \ {\rm and}\ y=0 ;\\
\{(u,v)|  u\perp x , \ v\perp y,  \ x_1\hat{u}+y_1v\in \mathbb{R} x \}   & \ {\rm if}\ x,y \in {\rm bd}{\cal K}\backslash\{0\} , x^Ty=0; \\
\{(u,v)| u\in \hat{y}^\circ,\ v\in \mathbb{R}_-\hat{y}  \} &  \ {\rm if}\ x=0, \ y\in {\rm bd}{\cal K}\backslash \{0\} ;\\
 \{(u,v)|  u\in \mathbb{R}_-\hat{x}, v\in \hat{x}^\circ   \}  & \ {\rm if}\ x\in {\rm bd}{\cal K}\backslash \{0\},\ y=0 ;\\
\{(u,v)|   u\in -{\cal K},\ v\in -{\cal K}  \}  & \ {\rm if}\ x=0, \ y=0.
\end{array}\right.
 \end{eqnarray*}
 \end{thm}

\noindent {\bf Proof.}
 We only  prove the case where  $x,y\in {\rm bd}{\cal K}\backslash \{0\}$ and $x^Ty=0$, since the other cases can be shown by using (\ref{regular-1}) and the expression for $\widehat{D}^*\Pi_{\cal K}$ given in \cite[Theorem 1]{OS08} by an elementary calculation.
In this case, by \cite[Lemma 2.3]{LZL14}, we have
 $x-y=\big((1-k)x_1, (1+k)x_2\big)$ with $k=y_1/x_1>0$. Note that $x-y\in (-{\cal K}\cup {\cal K})^c$. So
 according to \cite[Lemma 1(i)]{OS08},
 \[{\cal J} \Pi_{\cal K}(x-y)= \frac{1}{1+k}I+\frac{1}{2}\begin{bmatrix}
 -\frac{1-k}{1+k} & \bar{x}_2^T \\
 \bar{x}_2 & -\frac{1-k}{1+k}\bar{x}_2\bar{x}_2^T
 \end{bmatrix}.\]
 Hence by Proposition \ref{Prop2.2},
 \begin{eqnarray}
 (u,v)\in \widehat{N}_\Omega (x,y) \ \ &\Longleftrightarrow &\ \  -v\in \widehat{D}^*\Pi_{\cal K}(x-y)(-u-v)\nonumber\\
\ \ &\Longleftrightarrow &\ \
\left(\frac{1}{1+k}I+\frac{1}{2}\begin{bmatrix}
 -\frac{1-k}{1+k} & \bar{x}_2^T \\
 \bar{x}_2 & -\frac{1-k}{1+k}\bar{x}_2\bar{x}_2^T
 \end{bmatrix}\right)\left( \begin{array}{c}u_1+v_1 \\
 u_2+v_2 \end{array} \right) = \left( \begin{array}{c}v_1 \\ v_2 \end{array} \right) \nonumber\\
 &\Longleftrightarrow &
 \left\{ \begin{array}{l}
 u_1+\bar{x}_2^T (u_2+v_2) =v_1 \\
 \bigg[(1+k)(u_1+v_1)-(1-k)\bar{x}_2^T (u_2+v_2) \bigg]\bar{x}_2=2k v_2-2u_2.
  \end{array} \right. \label{equation}
 \end{eqnarray}
In what follows, we first show that the following inclusion holds
 \[
 \widehat{N}_\Omega (x,y) \subset \{(u,v)|\  v\perp y,\ u\perp x,\ x_1\hat{u}+y_1v\in \mathbb{R} x\}
 \] and then show the converse inclusion holds.
 Let $(u,v)\in \widehat{N}_\Omega (x,y)$.
  Take $x'\in {\rm bd}{\cal K}\backslash \{0\}$ and $y':=k\hat{x'}\in {\rm bd}{\cal K}\backslash \{0\}$. Then  $\langle x', y' \rangle =0$, i.e., $(x',y')\in \Omega$. Hence
 \begin{eqnarray}
 &&\frac{\langle (u,v), (x',y')-(x,y) \rangle }{\|(x',y')-(x,y)\|}=
 \frac{\langle u, x'-x \rangle +\langle v, y'-y \rangle  }{\|(x'-x,y'-y)\|}=\frac{\langle u, x'-x \rangle +\langle v, k\hat{x'}-k\hat{x} \rangle  }{\|(x'-x,k\hat{x'}-k\hat{x})\|} \nonumber\\
 &=& \frac{\langle u, x'-x \rangle +\langle k\hat{v}, x'-x \rangle  }{\|(x'-x,kx'-kx)\|}=\frac{\langle u, x'-x \rangle +\langle k\hat{v}, x'-x \rangle}{\sqrt{1+k^2}\|x'-x\|}=\frac{\langle u+k\hat{v}, x'-x \rangle}{\sqrt{1+k^2}\|x'-x\|}, \label{definition-1}
 \end{eqnarray}
 where we have used the fact that $\langle a, \hat{b} \rangle =\langle \hat{a} ,b \rangle  $ and $\|(a,\hat{b})\|=\|(a,b)\|$ for arbitrary vectors $a,b\in \mathbb{R}^m$.
 Since $(u,v)\in \widehat{N}_\Omega (x,y)$, it follows from (\ref{definition-1}) that
 \[
 \limsup\limits_{ {x'{\to }x}\atop {x'\in {\rm bd}{\cal K}\backslash \{0\}} } \frac{\langle u+k\hat{v}, x'-x \rangle}{\sqrt{1+k^2}\|x'-x\|}=\limsup\limits_{(x',y')\overset{\Omega}{\to}(x,y)}\frac{\langle (u,v), (x',y')-(x,y) \rangle }{\|(x',y')-(x,y)\|}
 \leq 0,
 \]
 which implies that
 \[
 u+k\hat{v}\in \widehat{N}_{{\rm bd}{\cal K}\backslash \{0\}}(x).
 \]
 Since $x_2\neq 0$,  ${\rm bd}{\cal K}=\{x| x_1-\|x_2\|=0\}$ is a smooth manifold near $x$. So $\widehat{N}_{{\rm bd}{\cal K}\backslash \{0\}}(x)=\{\mathbb{R} \hat{x}\}$ (see also \cite[Example 6.8]{rw}). Thus
$
 u+k\hat{v}\in \mathbb{R} \hat{x}.
$
 On the other hand, if in particular we choose $(x',y'):=(x,k'y)$ with $k'\to 1$. Then $(x',y')\in \Omega$ and
 \begin{eqnarray}\label{eq:1}
 &&\frac{\langle (u,v), (x',y')-(x,y) \rangle }{\|(x',y')-(x,y)\|}=
 \frac{\langle u, x'-x \rangle +\langle v, y'-y \rangle  }{\|(x'-x,y'-y)\|}=\frac{(k'-1)\langle v, y \rangle  }{|k'-1|\|y\|}. \end{eqnarray}
 Since $(u,v)\in \widehat{N}_\Omega (x,y)$,  it follows from the definition of regular normal cone and (\ref{eq:1}) that
 \[
 \limsup\limits_{k'\to 1}\frac{(k'-1)\langle v, y \rangle  }{|k'-1|\|y\|}=\limsup\limits_{(x',y')\overset{\Omega}{\to}(x,y)}\frac{\langle (u,v), (x',y')-(x,y) \rangle }{\|(x',y')-(x,y)\|}\leq 0,
 \]
 which implies that $v\perp y$. Similarly, we obtain $u\perp x$.

 From the above arguments, we have
 \begin{eqnarray*}
 \widehat{N}_\Omega (x,y) &\subset &\{(u,v)|\ u\perp x, \ v\perp y,\ u+k\hat{v}\in \mathbb{R} \hat{x}\}\\
 &=& \{(u,v)|\ u\perp x, \ v\perp y, \ x_1\hat{u}+y_1v\in \mathbb{R} x\}.
 \end{eqnarray*}
 Now we show that the converse inclusion holds. Let $(u,v)$ lie in the right hand side of the above inclusion. Then there exists $\beta \in \mathbb{R}$ such that
 \begin{eqnarray}
&&  \left(\begin{array}{c} u_1\\ u_2 \end{array}\right) + k\left(\begin{array}{c} v_1\\ -v_2 \end{array}\right)=\beta \left(\begin{array}{c} x_1\\ -x_2 \end{array}\right),\label{eq:2}
 \ \
 u_1 x_1+ u_2^T x_2=0, \ \ v_1y_1+v_2^T y_2=0,
 \end{eqnarray}
 which implies that $u_1+u_2^T\bar{x}_2=0$  and $v_1+v_2^T\bar{y}_2=0$
  since $x_1=\|x_2\|>0$ and $y_1=\|y_2\|>0$. Since   $\bar{x}_2=-\bar{y}_2$,  it follows that $v_1-v_2^T\bar{x}_2=0$. Hence $
 u_1+\bar{x}_2^T(u_2+v_2)=u_1+\bar{x}_2^T u_2+\bar{x}_2^T v_2= \bar{x}_2^T v_2=v_1
 $
 and
 \begin{eqnarray*}
 \bigg[(1+k)(u_1+v_1)-(1-k)\bar{x}_2^T (u_2+v_2) \bigg]\bar{x}_2
 &=&\bigg[(1+k)(u_1+v_1)-(1-k)(-u_1+v_1) \bigg]\bar{x}_2\\
 &=&\bigg[2u_1+2kv_1 \bigg]\bar{x}_2\\
 &=& 2\beta x_1 \bar{x}_2 = 2 \beta x_2= 2kv_2-2u_2,
 \end{eqnarray*}
 where the third and fifth equalities follow from (\ref{eq:2}). Thus
 $(u,v)$ satisfies (\ref{equation}), i.e., $(u,v) \in  \widehat{N}_\Omega (x,y)$.

 \BOX

\vskip 10 true pt

\section{Equivalence of the proximal and regular normal cones}
In this section we show that for the second-order cone complementarity set, the proximal normal cone coincides with the regular normal cone.
 Towards this end, we first show that the metric projection operator is not only B-differentiable but also calmly B-differentiable.
\begin{lemma}\label{calm} The metric projection operators $ \Pi_{\cal K}(\cdot)$ and $ \Pi_{{\cal K}^\circ}(\cdot)$
are  calmly B-differentiable for any given
$x\in \mathbb{R}^m$, i.e., for any $h\rightarrow 0,$
\begin{eqnarray*}
&&   \Pi_{\cal K}(x+h)-\Pi_{\cal K}(x)-\Pi'_{\cal K}(x;h)=O(\|h\|^2),\label{proximalP}\\
&&   \Pi_{{\cal K}^\circ}(x+h)-\Pi_{{\cal K}^\circ}(x)-\Pi'_{{\cal K}^\circ}(x;h)=O(\|h\|^2). \label{proximalPP}
  \end{eqnarray*}
 \end{lemma}

 \noindent {\bf Proof.}
We only prove the result  for $\Pi_{\cal K}$ since the proof for $\Pi_{{\cal K}^\circ}$ is exactly similar.  Consider the following six cases. \\
{\bf Case 1} $x\in {\rm int}{\cal K}$. In this case  $\Pi_{{\cal K}}(x)=x$, $\Pi_{{\cal K}}(x+h)=x+h$ for $h$ sufficiently close to $0$,  and $\Pi'_{{\cal K}}(x;h)=h$ by \cite[Lemma 2(i)]{OS08}. So
 \[
 \Pi_{{\cal K}}(x+h)-\Pi_{{\cal K}}(x)-\Pi'_{{\cal K}}(x;h)=x+h-x-h=0=O(\|h\|^2).
 \]
{\bf Case 2} $x\in -{\rm int}{\cal K}$. This case is symmetric to Case 1 and we omit the proof.

\noindent
 {\bf Case 3} $x\in  (-{\cal K}\cup {\cal K})^c$. Then  for $h$ sufficiently close to $0$, we have $x+h \in  (-{\cal K}\cup {\cal K})^c$ and so $\lambda_1(x)=x_1-\|x_2\|<0, \lambda_1(x+h)=(x_1+h_1)-\|x_2+h_2\|<0$.   By (\ref{proj}) and \cite[Lemma 2(i)]{OS08},
 \begin{eqnarray}
 \lefteqn{2\bigg[\Pi_{{\cal K}}(x+h)-\Pi_{{\cal K}}(x)-\Pi'_{{\cal K}}(x;h)\bigg]} \nonumber \\
 &=&(x_1+h_1+\|x_2+h_2\|)\left(\begin{array}{c}
 1\\ \frac{x_2+h_2}{\|x_2+h_2\|}\end{array}
 \right)-(x_1+\|x_2\|)\left(\begin{array}{c}
 1\\ \frac{x_2}{\|x_2\|}\end{array}\right) \nonumber\\
 &&-\begin{bmatrix}
 1 & \bar{x}_2^T \\
  \bar{x}_2 &  I+ \ds\frac{x_1}{\|x_2\|}\Big(I-\bar{x}_2\bar{x}_2^T\Big)
 \end{bmatrix} \left(\begin{array}{c}h_1 \\ h_2 \end{array}\right).\label{direction diff-1}
 \end{eqnarray}
The first component of the right hand side of  (\ref{direction diff-1}) is equal to
 \begin{eqnarray*}
\lefteqn{ x_1+h_1+\|x_2+h_2\|-(x_1+\|x_2\|)-(h_1+\bar{x}_2^T h_2)}\\
&&=
 \|x_2+h_2\|-\|x_2\|-\bar{x}_2^T h_2=O(\|h_2\|^2)=O(\|h\|^2),
 \end{eqnarray*}
 where the second equality holds by the fact that the norm is second-order continuously differentiable at $x_2\neq 0$.
 The second component of the right hand side of (\ref{direction diff-1}) is equal to
 \begin{eqnarray*}
 &&(x_1+h_1+\|x_2+h_2\|)\frac{x_2+h_2}{\|x_2+h_2\|}-(x_1+\|x_2\|)\frac{x_2}{\|x_2\|}-
 h_1\frac{x_2}{\|x_2\|}-h_2-\frac{x_1}{\|x_2\|}(I-\bar{x}_2 \bar{x}_2^T)h_2\\
 &=& x_1\left[ \frac{x_2+h_2}{\|x_2+h_2\|}-\frac{x_2}{\|x_2\|}-\frac{I-\bar{x}_2\bar{x}_2^T}{\|x_2\|}h_2 \right]+h_1\left[\frac{x_2+h_2}{\|x_2+h_2\|}-\frac{x_2}{\|x_2\|} \right]\\
 &=& x_1\left[ \frac{x_2+h_2}{\|x_2+h_2\|}-\frac{x_2}{\|x_2\|}-\frac{I-\bar{x}_2\bar{x}_2^T}{\|x_2\|}h_2 \right]+O(\|h\|^2)= O(\|h\|^2),
 \end{eqnarray*}
 where the second equality holds  by the Lipschitz continuity of $x_2/\|x_2\|$ and the last equality follows from the second-order continuous differentiability of $x_2/\|x_2\|$ at
 $x_2\neq 0$.

\noindent
{\bf Case 4} $x\in {\rm bd}{\cal K}\backslash\{0\}$. In this case $\lambda_1(x)=0$ and $\lambda_2(x) >0$. Then by (\ref{proj}) and \cite[Lemma 2(ii)]{OS08}, for $h$ sufficiently close to $0$,
\[
   \Pi_{{\cal K}}(x)=\frac{1}{2}\Big(x_1+\|x_2\|\Big)\left(\begin{array}{c}
 1\\ \frac{x_2}{\|x_2\|}\end{array}\right), \ \ \Pi'_{{\cal K}}(x;h)=h-\frac{1}{2}(h_1-\bar{x}_2^T h_2)_-\left(\begin{array}{c}1\\ -\bar{x}_2 \end{array}\right),
 \]
 and
 \[\Pi_{{\cal K}}(x+h)=\frac{1}{2}\Big(x_1+h_1-\|x_2+h_2\|\Big)_+\left(\begin{array}{c}
 1 \\ -\frac{x_2+h_2}{\|x_2+h_2\|}\end{array}\right)+\frac{1}{2}
 \Big(x_1+h_1+\|x_2+h_2\|\Big)
 \left(\begin{array}{c}
 1 \\ \frac{x_2+h_2}{\|x_2+h_2\|}\end{array}\right).\]
 Then the first component of $2\big[\Pi_{\cal K}(x+h)-\Pi_{{\cal K}}(x)-\Pi'_{\cal K}(x;h)\big]$ is
 \begin{eqnarray}
 &&\Big(x_1+h_1-\|x_2+h_2\|\Big)_+ +(x_1+h_1+\|x_2+h_2\|)-(x_1+\|x_2\|)-\Big(2h_1-(h_1-\bar{x}_2^T h_2)_-\Big)\nonumber\\
 &=& \Big( x_1+h_1-\|x_2+h_2\| \Big)- \Big(x_1+h_1-\|x_2+h_2\|\Big)_-+(x_1+h_1+\|x_2+h_2\|)-(x_1+\|x_2\|)\nonumber\\
 &&-\Big(2h_1-(h_1-\bar{x}_2^T h_2)_-\Big)\nonumber\\
 &=& - \Big(x_1+h_1-\|x_2+h_2\|\Big)_-+\Big(h_1-\bar{x}_2^T h_2\Big)_- \nonumber\\
 &=& - \Big(h_1-\bar{x}_2^Th_2+O(\|h_2\|^2)\Big)_-+\Big(h_1-\bar{x}_2^T h_2\Big)_- \nonumber\\
 &=& O(\|h_2\|^2)=O(\|h\|^2),\label{direction diff-2}
 \end{eqnarray}
 where the fourth equality holds since
 $ \big(h_1-\bar{x}_2^Th_2+O(\|h_2\|^2)\big)_-=\big(h_1-\bar{x}_2^T h_2\big)_-+O(\|h_2\|^2)
 $ by virtue of Lipschitz continuity of the function $t_-:=\min\{0,t\}$.
 According to (\ref{direction diff-2}) we have
 \begin{equation}\label{eq:calm-B-diff}
 -\Big(x_1+h_1-\|x_2+h_2\|\Big)_+=(x_1+h_1+\|x_2+h_2\|)-(x_1+\|x_2\|)
 -\Big(2h_1-(h_1-\bar{x}_2^T h_2)_-\Big)+O(\|h\|^2).
 \end{equation}
 The second component of $2\big[\Pi_{\cal K}(x+h)-\Pi_{{\cal K}}(x)-\Pi'_{\cal K}(x;h)\big]$ is
 \begin{eqnarray*}
 &&\bigg[-\Big(x_1+h_1-\|x_2+h_2\|\Big)_++\Big(x_1+h_1+\|x_2+h_2\|\Big)\bigg]\frac{x_2+h_2}{\|x_2+h_2\|}
 -(x_1+\|x_2\|)\frac{x_2}{\|x_2\|}-\\
 &&\ \ \Big(2h_2+(h_1-\bar{x}_2^Th_2)_-\bar{x}_2\Big)\\
 &=& \Big[ 2\|x_2+h_2\|+(h_1-\bar{x}_2^T h_2)_-+O(\|h\|^2)\Big]\frac{x_2+h_2}{\|x_2+h_2\|}
 -(x_1+\|x_2\|)\frac{x_2}{\|x_2\|}\\
 &&-\Big(2h_2+(h_1-\bar{x}_2^Th_2)_-\bar{x}_2\Big)\\
 &=& \Big(h_1-\bar{x}_2^T h_2\Big)_-\left[\frac{x_2+h_2}{\|x_2+h_2\|}-\frac{x_2}{\|x_2\|}\right]+O(\|h\|^2)
 = O(\|h\|^2),
 \end{eqnarray*}
 where the second equality follows from (\ref{eq:calm-B-diff}) and the last equality follows from the fact that $h_1-\bar{x}_2^Th_2=O(\|h\|)$ and the Lipschitz continuity of $x_2/\|x_2\|$ since $x_2\neq 0$ in this case.\\
{\bf Case 5} $x\in -{\rm bd}{\cal K}\backslash\{0\}$. In this case $\Pi_{\cal K}(x)=0$ and for $h$ that is very close to zero, $\lambda_1(x+h)<0$ and by (\ref{proj}) and \cite[Lemma 2(iii)]{OS08},
 \[\Pi_{\cal K}(x+h)=\frac{1}{2}\Big(x_1+h_1+\|x_2+h_2\|\Big)_+\left(\begin{array}{c}
 1\\ \frac{x_2+h_2}{\|x_2+h_2\|}\end{array} \right),\ \
 \Pi'_{\cal K}(x;h)=\frac{1}{2}(h_1+\bar{x}_2^Th_2)_+\left(\begin{array}{c}
 1 \\ \bar{x}_2
 \end{array}\right).
 \]
 The first component of $2\big[\Pi_{\cal K}(x+h)-\Pi_{{\cal K}}(x)-\Pi'_{\cal K}(x;h)\big]$ is
 \[
 \Big(x_1+h_1+\|x_2+h_2\|\Big)_+-\Big(h_1+\bar{x}_2^Th_2\Big)_+
 =
 \Big(h_1+\bar{x}_2^Th_2+O(\|h_2\|^2)\Big)_+-\Big(h_1+\bar{x}_2^Th_2\Big)_+=O(\|h_2\|^2).
 \]
 The second component of $2\big[\Pi_{\cal K}(x+h)-\Pi_{{\cal K}}x-\Pi'_{\cal K}(x;h)\big]$ is
 \begin{eqnarray*}
\lefteqn{\Big(x_1+h_1+\|x_2+h_2\|\Big)_+\frac{x_2+h_2}{\|x_2+h_2\|}-\Big(h_1+\bar{x}_2^Th_2\Big)_+\frac{x_2}{\|x_2\|}}\\
 &=&\Big(h_1+\bar{x}_2^T h_2\Big)_+\left(\frac{x_2+h_2}{\|x_2+h_2\|}-\frac{x_2}{\|x_2\|}\right)+O(\|h_2\|^2)= O(\|h\|^2),
 \end{eqnarray*}
 where the last equality follows from  $h_1+\bar{x}_2^Th_2=O(\|h\|)$ and the Lipschitz continuity of $x_2/\|x_2\|$ since $x_2\neq 0$ in this case.\\
 {\bf Case 6} $x=0$. Then $\Pi_{\cal \K}(x)=0$, $\Pi_{\cal \K}(x+h)=\Pi_{\cal \K}(h)$ and $\Pi'_{\cal \K}(x;h)=\Pi_{\cal \K}(h)$ by \cite[Lemma 2(iv)]{OS08}. Thus
 \[
 \Pi_{\cal K}(x+h)-\Pi_{{\cal K}}(x)-\Pi'_{\cal K}(x;h)=0=O(\|h\|^2).
 \]
 \BOX


  According to Lemma \ref{calm}, we can obtain the following result by using a similar proof technique as \cite[Proposition 3.1]{DSY14}.
 \begin{lemma}\label{prop3.2} Let $(x,y)\in \Omega$. Then
 $(u,v)\in N^\pi_\Omega(x,y)$ if and only if
 \begin{equation}\label{proximal normal}
 \langle u+v, \Pi'_{\cal K}(x-y;h) \rangle -\langle v,h\rangle\leq 0, \ \ \ \ \forall h\in \mathbb{R}^m.
 \end{equation}
 \end{lemma}

 With these preparations, the equivalence between the regular and proximal normal cone is given below.
\begin{thm}\label{equal}   Let $(x,y)\in \Omega:=\{(x,y)| x\in {\cal K}, y\in {\cal K}, x^Ty=0\}$. Then
 $\widehat{N}_\Omega(x,y)=N^\pi_\Omega(x,y)$.
 \end{thm}

 \noindent {\bf Proof.}
 Let $(x,y)\in \Omega$. Consider the following cases.

  \noindent {\bf Case 1} $x\in {\rm int}{\cal K}$,  $y=0$, or $x=0$, $y\in {\rm int}{\cal K}$, or $x,y\in {\rm bd}{\cal K}\backslash \{0\}$. In this case  $\Pi_{\cal K}$ is continuously differentiable at $x-y$.   By Lemma \ref{prop3.2}, $(u,v) \in N^\pi_\Omega(x,y)$ if and only if  (\ref{proximal normal}) holds. Since $\Pi_{\cal K}(x-y)$ is continuously differentiable at $x-y$,  (\ref{proximal normal}) takes the form
 \[
  \langle \nabla \Pi_{\cal K} (x-y)(u+v)-v, h \rangle \leq 0, \ \ \ \ \forall h\in \mathbb{R}^m,
 \]
 or equivalently,
 \[
 \nabla \Pi_{\cal K} (x-y)(u+v)-v=0.
 \]
 By Proposition \ref{Prop2.2},  the above equation holds if and only if  $(u,v) \in \widehat{N}_\Omega(x,y)$ and hence ${N}^\pi_\Omega(x,y)=\widehat{N}_\Omega(x,y)$.

 \noindent
{\bf Case 2} $x=0$ and $y\in {\rm bd}{\cal K}\backslash \{0\}$. In this case $x-y=-y\in -{\rm bd}{\cal K}\backslash \{0\}$. Hence by \cite[Lemma 2(iii)]{OS08} and the fact that $c_2(-y)=c_1(y)$, $\Pi'_{\cal K}(x-y;h)=2(c_1(y)^T h )_+ c_1(y)$. So (\ref{proximal normal}) takes the form
 \begin{eqnarray*}
  &&\left\langle u+v, 2(c_1(y)^T h )_+ c_1(y) \right\rangle -\langle v,h\rangle\leq 0 \ \ \ \ \forall h\in \mathbb{R}^m\\
  &\Longleftrightarrow &
  \left\{\begin{array}{cc}
  \langle -v,h \rangle  \leq 0 & \mbox{ if } \ c_1(y)^T h\leq 0 \\
  \langle u+v,2c_1(y)^T h c_1(y) \rangle -\langle v, h\rangle \leq 0 & \mbox{ if } \ c_1(y)^T h\geq 0
  \end{array} \right.\\
  &\Longleftrightarrow &
  \left\{\begin{array}{cc}
  \langle -v,h \rangle \leq 0& \ \mbox{ if } \ c_1(y)^T h\leq 0 \\
  \langle 2(u+v)^T c_1(y) c_1(y)-v, h\rangle \leq 0 & \ \mbox{ if } \ c_1(y)^T h\geq 0
  \end{array} \right.\\
  &\Longleftrightarrow &
  \exists \alpha,\beta \geq 0 \ \ such \ \ that \ \  -v=\alpha c_1(y)
  \ \ and \ \  2(u+v)^T c_1(y) c_1(y)-v=-\beta c_1(y) \\
  &\Longleftrightarrow &
  \exists \alpha,\beta \geq 0 \ \ such \ \ that \ \  -v=\alpha c_1(y)
  \ \ and \ \  2u^T c_1(y) c_1(y)=-\beta c_1(y) \\
  &\Longleftrightarrow &
  v\in \mathbb{R}_- c_1(y)
  \ \ and \ \  \langle u, c_1(y) \rangle  \leq 0.
  \end{eqnarray*}
  Since $y_1=\|y_2\|>0$, we have $c_1(y)=\frac{1}{2y_1}\hat{y}$ and hence $(u,v) \in {N}^\pi_\Omega(x,y)$ if and only if $u\in \hat{y}^\circ, v\in \mathbb{R}_-\hat{y}$. The equivalence of the two normal cones follows from  the exact formula of $\widehat{N}_\Omega(x,y)$ in Theorem \ref{formula-regular normal conenew}.

\noindent
{\bf Case 3} $x\in {\rm bd}{\cal K}\backslash \{0\}$ and $y=0$. In this case  $x-y=x$ and $c_1(x-y)=c_1(x)$. Hence by \cite[Lemma 2(ii)]{OS08}, $\Pi'_{\cal K}(x-y;h)=h-2\big(c_1(x)^Th\big)_- c_1(x)$. So (\ref{proximal normal}) takes the form
  \begin{eqnarray*}
  &&\left\langle u+v, h-2(c_1(x)^Th)_- c_1(x)\right\rangle -\langle v,h \rangle \leq 0, \ \ \forall h\in \mathbb{R}^m\\
  &\Longleftrightarrow &
  \left\{\begin{array}{cc}
  \langle u+v, h\rangle -\langle v,h \rangle\leq 0 &   \mbox{ if } \ c_1 (x)^T h\geq 0 \\
  \langle u+v,h-2c_1(x)^Th c_1(x) \rangle -\langle v, h\rangle  \leq 0 &  \mbox{ if } \ c_1(x)^T h\leq 0
  \end{array} \right.\\
  &\Longleftrightarrow &
  \left\{\begin{array}{cc}
  \langle u,h \rangle \leq 0&   \mbox{ if } \ c_1(x)^T h\geq 0 \\
  \langle u-2(u+v)^T c_1(x) c_1(x), h\rangle \leq 0 &  \mbox{ if } \ c_1(x)^T h\leq 0
  \end{array} \right.\\
  &\Longleftrightarrow &
  \exists \alpha,\beta \geq 0 \ \ such \ \ that \ \  u=-\alpha c_1(x)
  \ \ and \ \  u-2(u+v)^T c_1(x) c_1(x)=\beta c_1(x) \\
  &\Longleftrightarrow &
  \exists \alpha,\beta \geq 0 \ \ such \ \ that \ \  u=-\alpha c_1(x)
  \ \ and \ \  -2v^T c_1(x) c_1(x)=\beta c_1(x) \\
  &\Longleftrightarrow &
  u\in \mathbb{R}_- c_1(x)
  \ \ and \ \  \langle v, c_1(x) \rangle  \leq 0.
  \end{eqnarray*}
   Since $x_1=\|x_2\|>0$, we have $c_1(x)=\frac{1}{2x_1}\hat{x}$ and hence $(u,v) \in {N}^\pi_\Omega(x,y)$ if and only if $u\in \mathbb{R}_-\hat{x},v\in \hat{x}^\circ$. The equivalence of the two normal cones follows from  the exact formula of $\widehat{N}_\Omega(x,y)$ in Theorem \ref{formula-regular normal conenew}.

   \noindent
{\bf Case 4} $x=0$ and $y=0$.  In this case  (\ref{proximal normal}) takes the form
 \begin{eqnarray*}
  &&\langle u+v, \Pi_{{\cal K}}(h)\rangle -\langle v, h \rangle \leq 0, \ \ \forall h\in \mathbb{R}^m\\
  &\Longleftrightarrow & \langle u, \Pi_{{\cal K}}(h)\rangle - \langle v, \Pi_{{\cal K}^\circ}(h) \rangle \leq 0, \ \ \forall h\in \mathbb{R}^m\\
  &\Longleftrightarrow & u\in {\cal K}^\circ=-{\cal K} \ \ and \ \ v\in -{\cal K}.
  \end{eqnarray*}
The equivalence of the two normal cones follows from  the exact formula of $\widehat{N}_\Omega(x,y)$ in Theorem \ref{formula-regular normal conenew}. \BOX



\section{Expression of the limiting normal cone}

Due to the mistake in the formula for the regular normal cone when $x,y\in {\rm bd}{\cal  K}\backslash \{0\}$, the limiting normal cone given in \cite[Proposition 2.2]{LZL14}
also contains mistakes  for the cases where $x=0, y\in {\rm bd}{\cal K}\backslash \{0\}$, or   $x\in {\rm bd}{\cal K}\backslash \{0\}, y=0$, or $x=y=0$.
The formulas of the limiting normal cone given in \cite[Proposition 2.2]{LZL14} for these three cases are
 \begin{equation} N_\Omega(x,y)=\{(u,v)|u\in \mathbb{R}^m,\ v=0
  \ {\rm or}
  \   u \in (\mathbb{R}_+\hat{y})^\circ , \ v\in \mathbb{R}_- \hat{y} \}\ \ \ {\rm if}\ \ x=0, y\in {\rm bd}{\cal K}\backslash \{0\},\label{we2}\end{equation}
\begin{equation} N_\Omega(x,y)=\{(u,v)|  u=0, v\in \mathbb{R}^m \ {\rm or} \ u\in \mathbb{R}_- \hat{x},\
v \in (\mathbb{R}_+\hat{x})^\circ \} \ \ \ {\rm if}\ \ x\in {\rm bd}{\cal K}\backslash \{0\}, y=0,  \label{we3}
  \end{equation}
  and $ {\rm if} \ x=y=0$,
 \begin{eqnarray}
N_\Omega(x,y)&=&  \Big \{(u,v)| u\in -{\cal K},\ v\in -{\cal K} \ \ {\rm or} \  \ u\in \mathbb{R}^m,\ v=0 \ \
  {\rm or }\ \  u=0,\ v\in \mathbb{R}^m\nonumber \\
  &&  \qquad \ \ {\rm or}\ \  u\in \mathbb{R}_-\xi,\ v\in (\mathbb{R}_+\xi)^\circ  \ \mbox{ for some } \xi \in C \nonumber\\
  &&\qquad \ \ {\rm or}\ \ u\in (\mathbb{R}_+\xi)^\circ, \  v\in \mathbb{R}_-\xi \ \mbox{ for some } \xi \in C \nonumber\\
  && \qquad \  \mbox{ or }\   u \in \mathbb{R}\hat{\xi}, \ v\in \mathbb{R}\xi \ \mbox{ for some } \xi \in C \Big \} \ \  \label{we4}
 \end{eqnarray}
 where $C$ is defined as
 \begin{equation}
C:=\{(1,w)|\ w\in \mathbb{R}^{m-1},\ \|w\|=1\}.\label{C}
\end{equation}
  When $m=2$, it is easy to see that for any $\xi \in C$ and any $\alpha\in [0,1]$,
 $$u\perp \xi, v\perp \hat{\xi} \Longleftrightarrow u\in \mathbb{R}\hat{\xi}, v\in \mathbb{R}\xi \Longleftrightarrow u\perp \xi, v\perp \hat{\xi}, \alpha \hat{u}+(1-\alpha)v\in \mathbb{R}\xi.$$
 Hence when $m=2$, the limiting normal cone formula
 (\ref{we4})  at  $(x,y)=(0,0)$  is equivalent to our formula to be given in Theorem   \ref{formula-normal cone}. The following example illustrates that the formula (\ref{we2}) is not correct even when the dimension $m=2$ (similarly, (\ref{we3}) is not correct by symmetrical analysis) and the formula (\ref{we4}) is incorrect when $m$ is greater than $3$.
 \begin{example}
 1) For $x=(0,0), y=(1,1)\in bd{\cal K}\backslash \{0\}$, let $u=(1,1)$ and $v=(2,-2)$. Since $x-y=(-1,-1)\in -{\rm bd}{\cal K}\backslash \{0\}$, by \cite[Lemma 1(iii) and Theorem 2(iii)]{OS08}
 \[
 \frac{1}{2} \begin{bmatrix}
 1 & -1\\
 -1 & 1
 \end{bmatrix} \left( -u-v \right)
\in D^*\Pi_{\cal K}(x-y)(-u-v).
 \]
 Hence $-v \in D^*\Pi_{\cal K}(x-y)(-u-v).$ By Proposition \ref{Prop2.2}, it follows that $(u,v)\in N_\Omega(x,y)$. But $v\notin \mathbb{R}_- \hat{y}=\mathbb{R}_- (1,-1)$. Hence the formula (\ref{we2}) is incorrect. However $(u,v)$ satisfies the formula we proposed in Theorem  \ref{formula-normal cone} below, since $u=(1,1)\perp (1,-1)=\hat{y},\  v=(2,-2) =2\hat{y}\in \mathbb{R}\hat{y}.$

 2) For $x=y=(0,0,0)$, let $u=(1,-1,1)$ and $v=(0,0,1)$. Note that $u\notin \mathbb{R}(1,-w)$
 and $v\notin \mathbb{R}(1,w)$ with $\|w\|=1$, and hence $(u,v)$ does not belong to set proposed by the formula (\ref{we4}).
 However, by letting $\alpha=1/2$ and $w=(1,0)^T$, we have
 \[
 \frac{1}{2}\begin{bmatrix}
 1 & w^T\\
 w & I
 \end{bmatrix}(-u-v)=-v.
 \]
 Hence $-v\in D^*\Pi_{\cal K}(x-y)(-u-v)$, i.e., $(u,v)\in N_\Omega(x,y)$ by Proposition \ref{Prop2.2}. Take $\xi=(1,1,0)$. Then $u\perp \xi$, $v\perp \hat{\xi}$, and $\frac{1}{2}\hat{u}+\frac{1}{2}v=\frac{1}{2}(1,1,0)\in \mathbb{R} \xi$, i.e., $(u,v)$ satisfies the formula proposed in Theorem \ref{formula-normal cone} below.
 \end{example}

 \medskip

We now give a correct formula for the limiting normal cone of the second-order cone complementarity set. Note that the conditions $x_1\hat{u}+y_1v\in \mathbb{R} x$ and $\alpha \hat{u}+(1-\alpha)v\in \mathbb{R} \xi, \alpha\in [0,1]$ are redundant  when $m=2$.

 \begin{thm}\label{formula-normal cone}
   Let $(x,y)\in \Omega:=\{(x,y)|x\in {\cal K}, y\in {\cal K}, x^Ty=0\}$. Then
 \begin{eqnarray*}
 N_\Omega(x,y)= \widehat{N}_\Omega (x,y)=\left\{ \begin{array}{l}
 \{(u,v)| u\in \mathbb{R}^m, \ v=0 \}  \ {\rm if} \ x=0,\ y\in {\rm int}{\cal K}; \\
\{(u,v)| u=0, \ v\in \mathbb{R}^m \} \ {\rm if}\ x\in {\rm int}{\cal K},\ y=0 ;\\
 \{(u,v)| u\perp x,\ v\perp y,\  x_1\hat{u}+y_1v\in \mathbb{R} x\}  \ {\rm if} \  x, y\in {\rm bd}{\cal K}\backslash\{0\}.
 \end{array} \right.
 \end{eqnarray*}
 For $x=0, y\in {\rm bd}{\cal K}\backslash \{0\}$,
 $$ N_\Omega(x,y)=\{(u,v)|u\in \mathbb{R}^m,\ v=0
  \ \ {\rm or} \ \
 u\perp \hat{y},\ v\in \mathbb{R} \hat{y}
 \ \ {\rm or} \ \ \langle u, \hat{y} \rangle \leq 0, \ v\in \mathbb{R}_- \hat{y}\};$$
 for $x\in {\rm bd}{\cal K}\backslash \{0\}, y=0$,
 $$ N_\Omega(x,y)=\{(u,v)|  u=0, v\in \mathbb{R}^m \ \ {\rm or}\ \ u\in \mathbb{R} \hat{x},\
 v\perp \hat{x}\ \ {\rm or} \ \ u\in \mathbb{R}_-\hat{x},
  \langle v, \hat{x} \rangle \leq 0 \};$$
  for $x=y=0$,
 \begin{eqnarray*}
N_\Omega(x,y)&=&  \{(u,v)|\  u\in -{\cal K}, v\in -{\cal K} \ \mbox{\rm or }  u\in \mathbb{R}^m, v=0 \
 \mbox{\rm or } u=0, v\in \mathbb{R}^m\\
  &&  \  \mbox{\rm or } u\in \mathbb{R}_-\xi,\ v\in \xi^\circ \ \  \mbox{\rm or }\ u\in \xi^\circ, v\in \mathbb{R}_-\xi\\
   && \ \mbox{\rm or } u\perp \xi, \ v\perp \hat{\xi}, \
 \alpha \hat{u}+(1-\alpha)v\in \mathbb{R} \xi, \
  \alpha\in [0,1], \mbox{ for some } \xi \in C\}
 \end{eqnarray*}
 where  $C$ is defined as in (\ref{C}).
 \end{thm}

 \noindent {\bf Proof.}
 Consider the following cases.

  \noindent {\bf  Case 1} $x=0,y\in {\rm int}{\cal K}$, or $x\in {\rm int}{\cal K},y=0$ or $x,y\in {\rm bd}{\cal K}\backslash \{0\}$. In these cases, it is easy to prove since all points in $\Omega$ near  $(x,y)$
 belong to the same type and hence the regular normal cone and the limiting normal coincide.

\noindent {\bf  Case 2}  $x=0$ and $y\in {\rm bd}{\cal K}\backslash \{0\}$. Let $z:=x-y$. Then $z\in -{\rm bd}{\cal K}\backslash \{0\}$ and hence according to \cite[Theorem 2(iii)]{OS08},
 \begin{eqnarray*}
\lefteqn{ D^*\Pi_{\cal K}(z)(-u-v)}\\
  &=&\left\{O, \frac{1}{2} \begin{bmatrix}
  1 & \bar{z}_2^T\\
  \bar{z}_2 & \bar{z}_2\bar{z}_2^T
   \end{bmatrix}\right\}(-u-v)\cup \{z^*|\ z^*\in \mathbb{R}_+ c_2(z), \ \ \langle -u-v-z^*, c_2(z) \rangle\geq 0 \}.
 \end{eqnarray*}

Since $O\in D^*\Pi_{\cal K}(z)(-u-v)$,  it follows from (\ref{regular-2}) that
$(u,v)$ with $u\in \mathbb{R}^m$ and $v=0$ belongs to $N_\Omega(x,y)$.
 Take $\ds\frac{1}{2} \begin{bmatrix}
  1 & \bar{z}_2^T\\
  \bar{z}_2 & \bar{z}_2\bar{z}_2^T
   \end{bmatrix}(-u-v)\in D^*\Pi_{\cal K}(z)(-u-v)$. Since $\bar{z}_2=-\bar{y}_2$, the following equivalences hold.
 \begin{eqnarray*}
 &&-v=\frac{1}{2} \begin{bmatrix}
 1 & -\bar{y}^T_2\\
 -\bar{y}_2 & \bar{y}_2\bar{y}_2^T
 \end{bmatrix} \left( -u-v \right)
 \Longleftrightarrow  \left(\begin{array}{c} -v_1 \\ -v_2 \end{array} \right)=\frac{1}{2} \begin{bmatrix}
 1 & -\bar{y}^T_2\\
 -\bar{y}_2 & \bar{y}_2\bar{y}_2^T
 \end{bmatrix} \left( \begin{array}{c}-u_1 -v_1 \\
 -u_2 -v_2 \end{array} \right) \\
 &\Longleftrightarrow &
 \left\{ \begin{array}{l}
 u_1+v_1-(u_2+v_2)^T\bar{y}_2=2v_1\\
 -(u_1+v_1)\bar{y}_2+(u_2+v_2)^T\bar{y}_2\bar{y}_2=2v_2
  \end{array} \right. \\
  &\Longleftrightarrow& \left\{ \begin{array}{l}
 u_1-(u_2+v_2)^T\bar{y}_2=v_1\\
 -(u_1+v_1)\bar{y}_2+(u_1-v_1)\bar{y}_2=2v_2
  \end{array} \right.\\
  &\Longleftrightarrow& \left\{ \begin{array}{l}
 u_1-(u_2+v_2)^T\bar{y}_2=v_1\\
 -v_1\bar{y}_2=v_2
  \end{array} \right.\\
  &\Longleftrightarrow& \left\{ \begin{array}{l}
 u_1-u_2^T\bar{y}_2=0\\
 -v_1\bar{y}_2=v_2
  \end{array} \right.\\
   &\Longleftrightarrow& \left\{ \begin{array}{l}
 u\perp \hat{y}\\
 v\in \mathbb{R} \hat{y}.
  \end{array} \right.
 \end{eqnarray*}
It follows from  (\ref{regular-2}) that $\{(u,v)| u\perp \hat{y},v\in \mathbb{R} \hat{y}\}\subset N_\Omega(x,y)$.

For $ \{z^*|\ z^*\in \mathbb{R}_+ c_2(z), \ \ \langle -u-v-z^*, c_2(z) \rangle\geq 0 \}\subset D^*\Pi_{\cal K}(z)(-u-v)$ we have
 \begin{eqnarray*}
 \left\{\begin{array}{l}
 -v\in \mathbb{R}_+ c_2(z) \\
 \langle -u, c_2(z) \rangle \geq 0
 \end{array}  \right. \Longleftrightarrow
 \left\{\begin{array}{l}
 -v\in \mathbb{R}_+ \hat{y} \\
 \langle -u, \hat{y} \rangle \geq 0
 \end{array}  \right. \Longleftrightarrow  \left\{\begin{array}{l}
 v\in \mathbb{R}_- \hat{y} \\
 \langle u, \hat{y} \rangle \leq 0
 \end{array}  \right.
 \end{eqnarray*}
 where the second equivalence comes from the fact that $c_2(z)=\frac{1}{2}(1,-\bar{y}_2)=\frac{1}{2y_1}\hat{y}$ with $y_1>0$ since $y\in {\rm bd}{\cal K}\backslash \{0\}$.
It follows from  (\ref{regular-2}) that $\{(u,v)|  \langle u, \hat{y} \rangle \leq 0,\ v\in \mathbb{R}_- \hat{y}\}\subset N_\Omega(x,y)$.
 Combining the above possibilities,  we have
 \[N_\Omega(x,y)=\{(u,v)|\ u\in \mathbb{R}^m,\ v=0 \ \ {\rm or} \ \
 u\perp \hat{y},\ v\in \mathbb{R} \hat{y} \  \ {\rm or} \ \ \langle u, \hat{y} \rangle \leq 0,\ v\in \mathbb{R}_- \hat{y}\}.\]

 \noindent {\bf  Case 3} $x\in {\rm bd}{\cal K}\backslash \{0\}$ and $y=0$. The proof of this case is similar to Case 2.

 \noindent {\bf  Case 4} $(x,y)=(0,0)$. By \cite[Theorem 2(iv)]{OS08}, we have
 \begin{eqnarray*}
 D^*\Pi_{\cal K}(0)(-u-v)&=&\partial_B \Pi_{\cal K}(0)(-u-v)\cup \{z^*|\ z^*\in {\cal K},\ -u-v-z^*\in {\cal K}\}\\
 &&\cup \bigcup\limits_{\xi\in C} \{z^*|\ -u-v-z^*\in \mathbb{R}_+\xi, \ \langle z^*, \xi \rangle\geq 0 \}\\
 && \cup \bigcup_{\xi\in C} \{z^*|\ z^*\in \mathbb{R}_+\xi, \ \langle -u-v-z^*,\xi\rangle\geq 0 \}.
 \end{eqnarray*}

Since $O\in \partial_B \Pi_{\cal K}(0)$, it follows from (\ref{regular-2}) that $(u,v) $ with $v=0$ and $u\in \mathbb{R}^m$ belongs to $N_\Omega(x,y)$.

  Since $I\in \partial_B \Pi_{\cal K}(0)$, $(u,v)$ with  $u=0$ and $v\in \mathbb{R}^m$ belongs to $N_\Omega(x,y)$.

 Since $\alpha I+{\ds\frac{1}{2}}\begin{bmatrix}
 1-2\alpha & w^T \\
 w & (1-2\alpha)ww^T
 \end{bmatrix}\in \partial_B \Pi_{\cal K}(0)$
 for any $\alpha\in [0,1]$ and $\|w\|=1$, by virtue of  (\ref{regular-2}),
 \[
 \left(\alpha I+\frac{1}{2}\begin{bmatrix}
 1-2\alpha & w^T \\
 w & (1-2\alpha)ww^T
 \end{bmatrix}\right)(
 -u-v)=-v \Longrightarrow (u,v) \in N_\Omega(x,y),
 \]
 which can be rewritten equivalently as
 \begin{eqnarray}\label{(x,y)=(0,0)}
 \left\{\begin{array}{l}
 u_1+w^T(u_2+v_2)=v_1\\
 \alpha u_2+\alpha u_1 w=(1-\alpha)v_2-(1-\alpha)v_1 w
 \end{array}\right. \Longrightarrow (u,v) \in N_\Omega(x,y).
 \end{eqnarray}
 We now claim that the solution set of the system of two equations in (\ref{(x,y)=(0,0)}) is
 \begin{equation}\label{solution at (0,0)}
 \bigg\{(u,v)\left|\ u\perp \left(\begin{array}{c}1\\ w\end{array}\right),\ \ v\perp \left(\begin{array}{c}1\\ -w\end{array}\right), \ \
 \alpha \hat{u}+(1-\alpha)v\in \mathbb{R} \left(\begin{array}{c}1\\ w\end{array}\right)\right.\bigg\}.
 \end{equation}
 Multiplying $w$ to the second equation in  the system of two equations in (\ref{(x,y)=(0,0)}) yields \begin{eqnarray*}
 v_2^Tw-v_1+\alpha v_1 =\alpha (u_2+v_2)^T w+\alpha u_1
 =\alpha (v_1-u_1)+\alpha u_1 = \alpha v_1
 \end{eqnarray*}
 where the second equality holds by the first equality  in (\ref{(x,y)=(0,0)}). This means that $v_2^Tw-v_1=0$, i.e., $v\perp (1,-w)$. Applying this  to the first equation in (\ref{(x,y)=(0,0)}) yields $u_1+w^T u_2=0$, i.e., $u\perp (1,w)$.
 Using (\ref{(x,y)=(0,0)}) again yields $
 (1-\alpha)v_2-\alpha u_2=\Big[ \alpha u_1+(1-\alpha)v_1\Big]w.$
 Let $\eta:=\alpha u_1+(1-\alpha)v_1$. Then
 $\alpha \hat{u}+ (1-\alpha)v=\eta(1,w)\in \mathbb{R} (1,w)$.
 Conversely, take $(u,v)$ satisfying (\ref{solution at (0,0)}), i.e., there exists $\eta\in \mathbb{R}$ such that
 \[
 \alpha \left(\begin{array}{c} u_1\\ -u_2 \end{array}\right)+
 (1-\alpha)\left(\begin{array}{c} v_1\\ v_2 \end{array}\right)=
 \eta \left(\begin{array}{c} 1\\ w \end{array}\right), \ \
  u\perp \left(\begin{array}{c}1\\ w\end{array}\right),\ \ v\perp \left(\begin{array}{c}1\\ -w\end{array}\right).\]
 Then
 \[
 u_1+(u_2+v_2)^Tw= u_1+ u_2^T w+v_2^T w=v_2^Tw=v_1
 \]
 and
 \begin{eqnarray*}
 \alpha u_2+\alpha u_1 w &=& -\eta w+ (1-\alpha)v_2+\Big[ \eta-(1-\alpha)v_1\Big]w\\
 &=&(1-\alpha)v_2-(1-\alpha)v_1w,
 \end{eqnarray*}
i.e.,  $(u,v)$ satisfies the system of equations in  (\ref{(x,y)=(0,0)}). It follows that  any element $(u,v)$ in the set (\ref{solution at (0,0)}) belongs to the limiting normal cone $N_\Omega(x,y)$.

 Since $\{z^*|\ z^*\in {\cal K},\ -u-v-z^*\in {\cal K}\}\subset D^*\Pi_{\cal K}(0)(-u-v)$,  by (\ref{regular-2}) any $(u,v)$ such that $v\in -{\cal K}$ and $u\in -{\cal K}$ lies in $N_\Omega(x,y)$. Similarly, from  $\{z^*|\ -u-v-z^*\in \mathbb{R}_+ \xi, \ \langle z^*,\xi \rangle\geq 0 \}$ we derive  that any $(u,v)$ such that $u\in \mathbb{R}_-\xi$ and $v\in \xi^\circ$ lies in $N_\Omega(x,y)$  and from $\{z^*|\ z^*\in \mathbb{R}_+\xi, \langle -u-v-z^*, \xi\rangle\geq 0 \}$ we derive  that any $(u,v)$ such that $v\in \mathbb{R}_-\xi$ and $u\in \xi^\circ$  lies in $N_\Omega(x,y)$.
Combining all possibilities yields the formula of $N_\Omega(x,y)$ at $(0,0)$. \BOX

\section*{Acknowledgements.} The authors are indebted to the two
anonymous reviewers for their  useful comments which helped  us to  make the paper more concise.



\end{document}